\documentclass[12pt]{article}

\usepackage{amssymb,amsmath,graphics}

\numberwithin{equation}{section}

\allowdisplaybreaks

\newenvironment{Proof}{\removelastskip\par\medskip
\noindent{\em Proof.}
\rm}{\penalty-20\null\hfill$\square$\par\medbreak}

\def\real{{\mathord{\mathbb R}}}

\newtheorem{prop}{Proposition}

\def\Dom{{\mathrm{{\rm Dom ~}}}}

\def\deltaa{{\mathord{\mathrm {Delta}}}}
\def\rhoa{{\mathord{{\mathrm {Rho}}}}}
\def\vegaa{{\mathord{{\mathrm {Vega}}}}}
\def\gammaa{{\mathord{{\mathrm {Gamma}}}}}
\def\thetaa{{\mathord{{\mathrm {Theta}}}}}

\begin{document}
 ~
\begin{center}
{\Large Price sensitivities for a general stochastic volatility model}
\\
		\bigskip
		{\large Youssef El-Khatib}\footnote{UAE University, Department of Mathematical Sciences, Al-Ain, P.O. Box 15551. United Arab Emirates.{\ E-mail : Youssef\_Elkhatib@uaeu.ac.ae.}
		}\ \ \ \ \ {\large Abdulnasser Hatemi-J}\footnote{UAE University, Department of Economics and Finance, Al-Ain, P.O. Box 15551, United Arab Emirates. {\ E-mail : Ahatemi@uaeu.ac.ae.}}
	\end{center}
~
\begin{abstract}
\baselineskip0.6cm
We deal with the calculation of price sensitivities for stochastic volatility models.  General forms for the dynamics of the underlying asset price and its volatility are considered. We make use of the chaotic (or Malliavin)
calculus to compute the price sensitivities. The obtained results are applied to several recent stochastic volatility models as well as the existing ones that are commonly used by practitioners. Each price sensitivity is a source of financial risk. The suggested formulas are expected to improve on the hedging of the underlying risk. 
\end{abstract}

%\noindent{{\em JEL Classification:}} G12, G10, C60.

\baselineskip0.77cm
\section{Introduction}
Mathematical tools are increasingly utilized by investors and financial institutions
in order to neutralize or reduce the underlying financial risk. Among others, price
sensitivities are commonly used in markets for financial derivatives in order to hedge against risk. This is indeed an active field of research. Recently, several papers have dealt with this important issue. It is shown in the literature that the valuation of financial derivatives is more accurate if the underlying data generating process is characterized by a stochastic  volatility process. The existing literature suggests using a stochastic process for the volatility in order to determine the price of financial derivatives. The aim of this paper is to deal with the derivation of the price sensitivities.  We derive price sensitivities for a general stochastic volatility model by making use of the Malliavin calculus. It should be mentioned that each financial trading position that is based on financial instruments has five price sensitivities, which are known as the Greeks in the financial literature. The precise calculation of these price sensitivities is of paramount importance for the immunization of potential financial risk of a financial trading position or a portfolio. The first price sensitivity is known as Delta, which is equal to the change of the trading position with respect to the price of the underlying asset under the ceteris paribus condition. The second price sensitivity is called Gamma, which is representing the change of the Delta for a portfolio of options with respect to the price of the underlying asset. The third source of risk in this situation is known as Vega and it embodies the change of the trading position with regard to a marginal change in the volatility of the original asset. The fourth price sensitivity is called Theta, which captures the change of the value of the underlying portfolio with respect to the time factor. The last price sensitivity is commonly known as Rho in the literature and it is capturing the sensitivity of the trading position with respect to the risk free rate (i.e. the interest rate). Each of these Greek measures represent a source of risk for the portfolio and traders must calculate the pertinent Greeks for their portfolio at the end of every trading day so as to take required action if the internal risk is higher than the pre-determined levels by the underlying financial institution.
The Malliaivn calculus is especially useful in this case because the price of a financial derivative is regarded as a stochastic process that does not have a closed form solution in general. Consequently, dealing with the price sensitivities via this method is an appropriate approach. Through the Malliavin calculus we are able to transform the differentiation into integration and thus provide an unbiased measure of each price sensitivity. Many of the existing contributions to the price sensitivities are based on the finite difference approach, which can indeed be considered an biased methodology. Conversely, the Malliavin method is unbiased and it can also be less time-consuming in terms of convergence.
The idea to make use of the Malliavin calculus for computing price sensitivities originates from \cite{fournie1999}. This first application was within the context of a market that is characterized by information generated by the Brownian motion.  Their method relies on the Malliavin derivative on the Wiener space, which contains two parts. The first part is the application of the chain rule and the second part is utilizing the fact that the derivative has an adjoint (Skorohod integral) which can be described by the Ito integral for adapted processes. Lately, this approach has been utilized by \cite{elkhatibhatemi11} for markets under stress or experiencing a financial crisis. There are also several other papers that use the Malliavin calculus for markets with jumps\footnote{We can find recent works using numerical techniques for pricing and hedging financial derivatives in jump markets see for instance \cite{elkhatibmdalal12} and \cite{elkhatibhajji13}.}. For example, \cite{elkhatibprivault04} uses the Poisson noise via the jump times. \cite{BBM07} takes into account both the jump timing as well as the amplitude of the underlying jumps.  \cite{DJ06} deals with jump-diffusion models via the Malliavin calculus with regard to the Brownian motion conditional on the Poisson component. In addition, \cite{BM06} permits the Poisson noise to take into account the amplitude of the jumps. However, the timing of the jumps is not taken into account. Conversely, \cite{elkhatibhatemi12} provide the price sensitivities by utilizing the Malliavin derivative on the Wiener space and the Poisson noise when the timing of the jumps is taken into account. This paper provides price sensitivities for a general stochastic volatility model that encompasses a number of well-known existing models as well as several new ones.$\\$  $\\$
After this introduction the remaining part of the paper is organized as follows:
Section two presents the model and an introduction to Malliavin calculus. Section
three deals with deriving the price sensitivities via the Malliavin approach. A general formula is
provided that encompasses different stochastic volatility models. The formula is applied to find the price sensitivities for several new stochastic volatility models as well as for other specific models that are well-known in the
literature. In the last Section concluding remarks are provided.
\section{Preliminaries}
In this section, we describe the general stochastic volatility model and we present some tools from
Malliaivin calculus needed to our study. $\newline$
\noindent
We consider two independent Brownian motion $(B_t)_{t\in [0,T]}$ and $(B^{'}_t)_{t\in
[0,T]}$ and a filtered
probability space $(\Omega, {\mathcal{F}},({\mathcal{F}}_t)_{t\in
	[0,T]}, P)$, where $({\mathcal{F}}_t)_{t\in [0,T]}$ is the natural filtration generated by $B$ and $B^{'}$.
Next we introduce a general framework for the stochastic volatility
model.
\subsection{The stochastic volatility model}
We assume that the marketplace contains only two assets. The first asset is a risk
free asset given by $(A_t:=e^{\int_0^t r_s ds})_{t\in [0,T]}$ where $r$ is the interest rate.
The second asset is the underlying risky asset on which a European call option is built.
It is assumed to have a stochastic volatility i.e the volatility is driven by a stochastic volatility model. The price of the underlying asset $(X_t)_{t\in [0,T]}$ and the volatility process are stochastic processes driven by the following stochastic differential equations
\begin{eqnarray}
\label{asset}
dX_t&=&\mu_t X_t dt+\sigma(Y_t) X_t dB_t,\\
\label{volatility}
dY_t&=&g(Y_t)dt+\beta[\rho dB_t+\sqrt{1-\rho^2}dB^{'}_t ],\     \   	 \ 	    \  t \in[0,T],
\end{eqnarray}
with $X_0=x>0$ and $Y_0=y \in \real$. $\mu$ is a deterministic function, $\beta \in \real$, $\rho \in [-1,1]$ and $\sigma \in
{\mathcal C}^2([0,T]\times \real)$ such that for any $t\in [0,T]$,
$\sigma(.)\neq 0$. $\sigma$ is the volatility of the underlying asset, $\beta$ measures the volatility of the volatility and $\rho$ represents a measure of dependency between the price of the underlying asset and its volatility. The market considered here is incomplete. There
is an infinity of E.M.M -Equivalent Martingale Measure- (i.e a
probability equivalent to $P$ under which the actualized price $(X_t
e^{-rt})_{t \in [0,T]}$ is a martingale). Let $Q$ be a fixed
$P$-E.M.M. $Q$ is identified by its Radon-Nikodym density w.r.t $P$, denoted
$\zeta_T$ and given by
$$ \zeta_T=\exp\left(\int_0^T a_s dB_s+b_s dB^{'}_s
-\frac{1}{2}\int_0^T (a^2_s+b^2_s)ds\right),$$ where
$(a_{t})_{t\in [0,T]}$ and $(b_{t})_{t\in [0,T]}$ are two
predictable processes s.t. 
\begin{equation}
\label{at}
a_t= -\frac{\mu_t-r_t}{\sigma(Y_t)},
\end{equation} and $b$
is arbitrary. Now let, for any $t\in
[0,T]$, $W_t=B_t -\int_0^t a_s ds$ and
$W^{'}_t=B^{'}_t -\int_0^t b_s ds$ then by the Girsanov theorem
$W$ and $W^{'}$ are two $Q$-Brownian motions. In the following we
work with a fixed $P$-E.M.M. $Q$ and we will use $E[.]$ (instead of
$E_{Q}$[.]) as the expectation under the probability $Q$. We have,
under $Q$, for any $t\in [0,T]$
\begin{eqnarray}
\label{Stock1}
dX_t&=&r_t X_t dt+\sigma(Y_t) X_t dW_t,\\
\label{Stock2}
dY_t&=&h(Y_t, X_t)dt+\beta[\rho dW_t+\sqrt{1-\rho^2}dW^{'}_t ],
\end{eqnarray}
where 
\begin{equation}
\label{hyt}
h(Y_t, X_t):=g(Y_t)+\beta \rho a_t +\beta \sqrt{1-\rho^2} b_t,
\end{equation}
$a_t$ is given by (\ref{at}) and $b_t$ is a predictable process, assumed depending on $Y_t$ and $X_t$ i.e. $b_t=b_t(Y_t, X_t)$.
\subsection{Malliavin derivative}
We give an introduction to Malliavin derivative
in Wiener space and we list some important results \footnote{The reader can refer to \cite{oksendal1996} for more detailed description on Malliaivin calculus.}.
We denote by $(D_t)_{t\in [0,T]}$ the Malliavin derivative on the direction of $W$. We denote by
$\mathbb{P}$ the set of random variables $F:\Omega \rightarrow
\real$, such that F has the representation
$$F(\omega)=f\left(\int_0^T f_1(t)dW_t,
\ldots,\int_0^T f_n(t)dW_t\right),$$ where
$f(x_1,\ldots,x_n)=\sum_{\alpha}a_{\alpha}x^{\alpha}$ is
a  polynomial in $n$ variables $x_1,\ldots,x_n$ and
deterministic functions  $f_i \in L^2([0,T])$ .
Let $\|.\|_{1,2}$ be the norm
$$\|F\|_{1,2}:=\|F\|_{L^2(\Omega)}+\|D_{\cdot} F\|_{L^2([0,T]\times\Omega)}, \ \ \ F \in L^2(\Omega).$$
 Then $\Dom(D)$, the domain of $D$, is equal to
$\mathbb{P}$ w.r.t the norm $\|.\|_{1,2}$. The next propositions are very useful when using the Malliavin derivative.
\begin{prop}
	\label{p1} Given $F=f\left(\int_0^T f_1(t)dW_t, \ldots,\int_0^T
	f_n(t)dW_t\right)\in \mathbb{P}$. We have
	$$
	D^W_t F =\sum_{k=0}^{k=n}\frac{\partial f}{\partial
		x_k}\left(\int_0^T f_1(t)dW_t,\ldots,\int_0^T
	f_n(t)dW_t\right)f_k(t).
	$$
\end{prop}
To calculate the Mallaivin derivative for It\^o integral, we will
use the following Proposition.
\begin{prop} \label{derivint} Let $(u_t)_{t\in [0,T]}$ be a
	${\mathcal{F}}_t-$adapted process, such that $u_t \in \Dom(D)$,
	we have
	$$D_t \int_0^T u_s dW_s=\int_t^T (D_t u_s)dW_s+ u_t.$$
\end{prop}
From now on, for any stochastic process $u$ and for $F\in \Dom(D)$
such that $u_. D_. F \in L^2([0,T])$  we let
\begin{equation}
\label{Du}
D_u F:=\langle D F,u\rangle_{L^2([0,T])}:=\int_0^T u_t D_t
F dt.
\end{equation}
 The next proposition presents some  important results that link $D$ and its adjoint $\delta$, known as the Skorohod integral. 
\begin{prop}
	a) Let $u \in \Dom(\delta)$ and $F\in \Dom(D)$, we have $E[D_u
	F]\leq C(u)\|F\|_{1,2}$, and $E[F\delta(u)]=E[D_u F]$.\\
	b) Consider a $L^2(\Omega \times [0,T])$-adapted stochastic process
	$u=(u_t)_{t\in [0,T]}$. We have $\delta(u)=\int_0^T u_t dW_t.$\\
	c) Let $F\in \Dom(D)$ and $u\in \Dom(\delta)$ such that $uF \in
	\Dom(\delta)$ thus $\delta(uF)=F\delta(u)-D_u F.$
\end{prop}
\section{Price sensitivities}
We consider a European option with payoff $f(X_T)$,
where $(X_t)_{t\in [0,T]}$ denotes the underlying asset price given
by the general stochastic volatility model (\ref{Stock1}-\ref{Stock2}). We denote by $C$ the value of the European option. We will compute the following price sensitivities:
\begin{eqnarray*}
	& & \deltaa = \frac{\partial C}{\partial x},
	\quad \quad
	\gammaa =\frac{\partial^2 C}{\partial x^2},
	\quad \quad \rhoa=\frac{\partial C}{\partial r},
	\quad \quad
	\vegaa = \frac{\partial C}{\partial \sigma},
\end{eqnarray*}
The last price sensitivities $\thetaa = \frac{\partial V}{\partial t}$ can be obtained using the partial differential equation satisfied by $C$, the price of the option. In the following proposition we find the Malliaivn derivatives of $X_T$ and $Y_t$ w.r.t $D$.
\begin{prop}
	\label{der01}
	For $0\leq t\leq T$, We have 
	\begin{equation}
\label{DX}
	D_t X_T=X_T \left(\sigma(Y_t)-\int_t^T \sigma^{'}(Y_s) \sigma(Y_s) D_t Y_s ds +\int_t^T \sigma^{'}(Y_s) D_t Y_s d{W}_s\right)
	\end{equation}
	and
		\begin{equation}
		\label{DY}
	D_t Y_s=\beta\rho+\int_t^s D_t h(Y_\nu, X_\nu)d\nu,
		\end{equation}
where
\begin{equation}
\label{hyx}
	D_t h(Y_\nu, X_\nu)=g^{'}(Y_\nu)D_t Y_\nu+\beta \rho (\mu_\nu-r_\nu)\frac{\sigma^{'}(Y_\nu)}{\sigma^2(Y_\nu)}D_t Y_\nu+\beta \sqrt{1-\rho^2} D_t b_\nu,
\end{equation}
with $0\leq t \leq \nu \leq s \leq T$.
\end{prop}
\begin{Proof}
The equality (\ref{DX}) can be obtained by applying the Malliavin derivative to (\ref{Stock1}). Then we use the chain rule, proposition~ 2, and 
\begin{equation}
\label{DdetInt}
D_t \int_0^T u_s ds=\int_t^T (D_t u_s)ds,
\end{equation}
when $(u_s)_{s\in [0,T]}$ is an adapted process. $\\$
To find $D_t Y_v$, we have from (\ref{Stock2}), for $0\leq t\leq v\leq T$,
\begin{eqnarray*}
	D_t Y_v &=&D_t\left(Y_0+\int_0^v h(Y_s, X_s)ds+\beta\rho W_v+\beta\sqrt{1-\rho^2}W^{'}_v\right)
	\\
	&=&\beta\rho+\int_t^v D_t h(Y_s, X_s)ds,
\end{eqnarray*}
where 
\begin{eqnarray*}
D_t h(Y_s, X_s)&=& D_t(g(Y_s)+\beta \rho a_s +\beta \sqrt{1-\rho^2} b_s)\\
&=&g^{'}(Y_s)D_t Y_s+\beta \rho D_t \frac{r_s-\mu_s}{\sigma(Y_s)} +\beta \sqrt{1-\rho^2} D_t b_s\\
&=& g^{'}(Y_s)D_t Y_s+\beta \rho (\mu_s-r_s) \frac{\sigma^{'}(Y_s)}{\sigma^2(Y_s)}D_t Y_s+\beta \sqrt{1-\rho^2} D_t b_s.
\end{eqnarray*}
\end{Proof}
The second and the third order derivatives of $X_T$ w.r.t $D$, essential for computing the different price sensitivities are given in the following proposition. 
\begin{prop}
	\label{der02}
	For $0\leq t\leq T$, we let $L_t^T:=\sigma(Y_t)-\int_t^T \sigma^{'}(Y_\nu) \sigma(Y_\nu) D_t Y_\nu d\nu +\int_t^T \sigma^{'}(Y_\nu) D_t Y\nu d{W}_\nu$, then we have $D_t X_T=X_T L_t^T$ and
	\begin{eqnarray}
	\label{der1} D_u X_T&=&X_T \int_0^T u_t L_t^Tdt\\
	\label{der2} D_u D_u X_T&=&X_T\left(\left(\int_0^T u_t L_t^T dt\right)^2+\int_0^T\int_s^Tu_s u_t D_sL_t^Tdtds\right)\\
	\nonumber D_u D_u D_u X_T &=& X_T\left(\left(\int_0^T u_t L_t^Tdt\right)^3+3\int_0^T u_t L_t^Tdt \int_0^T\int_s^Tu_su_tD_s L_t^Tdtds\right.\\
	&&
	\label{der3}
	\left.+\int_0^T\int_r^T\int_s^T u_ru_su_t D_r D_s L_t^T dtdsdr\right)
	\end{eqnarray}
\end{prop}
\begin{Proof}
The equalities (\ref{der1}-\ref{der3}) are obtained using (\ref{Du})  and the chain rule of the Malliavin derivative. 
\end{Proof}
The derivatives $D_s L_t^T$, $D_r D_s L_t^T$, $D_s D_t Y_v$ and $D_r D_s D_t Y_v$ can be found using the chain rules of the Mallaivin derivative. Similar calculations can be found in \cite{elkhatib2009}.
\subsection{First order price sensitivities: $\deltaa$, $\rhoa$, $\vegaa$}
Let $C=E[f(X_T^{\zeta})]$ be the price of the option, $\zeta$ can take values: the asset price $x$ to obtain $\deltaa$, the interest rate $r$  for $\rhoa$, and $\sigma$ for $\vegaa$. We have 
$$
\frac{\partial}{\partial \zeta}E\left[
f(X_T^\zeta)
\right]
= E\left[f(X_T^\zeta)\left(\frac{\partial_\zeta X_T^\zeta}
{D_u X_T^\zeta }\delta(u)-D_u
\left(\frac{X_T^\zeta \partial_\zeta X_T^\zeta}{D_u X_T^\zeta
}\right)\right)\right].
$$
Next we compute the $\deltaa$, the $\rhoa$ and $\vegaa$ can be computed by the
same way. The $\deltaa$ corresponds to $\zeta=x$, so
$\partial_\zeta S_T=\partial_x S_T=\frac{1}{x}S_T$ and we have
\begin{eqnarray*}
	\nonumber
	\deltaa&=&E\left[f(X_T)\left(\frac{\partial_x X_T}
	{D_u X_T}\delta(u)-D_u
	\left(\frac{X_T\partial_x X_T}{D_u X_T
	}\right)\right)\right]\\
	\nonumber
	&=&E\left[f(X_T)\left(\frac{X_T}
	{x D_u X_T}\delta(u)-D_u
	\left(\frac{X_T^2}{xD_u X_T
	}\right)\right)\right]\\
	\nonumber
	&=&\frac{1}{x}E\left[f(X_T)\left(\frac{1}
	{\int_0^T u_t L_t^Tdt}\delta(u)-2X_T+
	\frac{X_T^2D_u D_u X_T}{(D_u X_T)^2
	}\right)\right]\\
	\nonumber
	&=&\frac{1}{x}E\left[f(X_T)\left(\frac{1}
	{\int_0^T u_t L_t^T dt}\delta(u)-2X_T\right.\right.\\
	&&\left.\left.+
	\frac{X_T(\int_0^T u_t L_t^T dt)^2+X_T\int_0^T u_s D_s(\int_0^T u_t L_t^T dt)ds}{(\int_0^T u_t L_t^Tdt)^2
	}\right)\right]\\
	\label{del}
	&=&\frac{1}{x}E\left[f(X_T)\left(\frac{1}
	{\int_0^T u_t L_t^T dt}\delta(u)-X_T\left(1-
	\frac{\int_0^T u_s (\int_s^T u_t D_s L_t^Tdt)ds}{(\int_0^T u_t L_t^Tdt)^2
	}\right)\right)\right].
\end{eqnarray*}
\subsection{$\gammaa$}
The $\gammaa$ is computed using the second order derivative of
$C=E[f(S_T)]$ w.r.t $x$ given by
\begin{eqnarray*}
	\label{gam}
	\lefteqn{\frac{\partial^2}{\partial x^2}
		E\left[f(X_T) \right]=\frac{\partial}{\partial x}\deltaa=\frac{1}{x}\frac{\partial}{\partial x}E\left[f(X_T)H\right]}\\
	\nonumber
	&=&\frac{1}{x}E\left[f(X_T)\left(\frac{H\partial_x X_T}
	{D_u X_T}\delta(u)-D_u\left(\frac{H\partial_x X_T}{D_u X_T}\right)+\partial_x H \right)\right]\\
	\nonumber
	&=&\frac{1}{x}E\left[f(X_T)\left(\frac{H X_T}
	{xD_u X_T}\delta(u)-\left(\frac{H(D_u X_T)^2+X_TD_uX_TD_uH+HX_TD_u D_u X_T}{x(D_u X_T)^2}\right)\right.\right.\\
	&&\left.\left.+\partial_x H \right)\right],
\end{eqnarray*}
where
\begin{eqnarray*}
	H&=&\frac{1}{\int_0^T u_t L_t^Tdt}\delta^W(u)-X_T\left(1-
	\frac{\int_0^T u_s (\int_s^T u_t D_s L_t^Tdt)ds}{(\int_0^T u_t L_t^T dt)^2
	}\right)\\
	D_u H&=&-\frac{\int_0^T u_s (\int_s^T u_t D_s L_t^Tdt)ds}{(\int_0^T u_t L_t^T dt)^2
	}\delta(u)\\
	&&-X_T\left(\int_0^T u_t L_t^T dt-\frac{\int_0^T u_s (\int_s^T u_t D_s L_t^Tdt)ds}{\int_0^T u_t L_t^T dt
	}\right.\\
	&&\left.-\frac{\int_0^T u_r\int_r^T u_s (\int_s^T u_t D_rD_s L_t^Tdt)ds}{(\int_0^T u_t L_t^T dt)^2}
	+2\frac{\left(\int_0^T u_s (\int_s^T u_t D_s L_t^T dt)ds\right)^2}{(\int_0^T u_t L_t^T dt)^3}\right)\\
	\partial_x H &=&-\frac{1}{x}X_T\left(1-\frac{\int_0^T u_s (\int_s^T u_t D_s L_t^T dt)ds}{(\int_0^T u_t L_t^T L_t^T dt)^2
	}\right).
\end{eqnarray*}
The previous calculations for price sensitivities are valid for most of the recent and well-known stochastic volatility models existing in the literature.  
\subsection{Applications}
We start first by the linear stochastic volatility model of \cite{jak13}. The underlying asset price $X_t$ is assumed to have the dynamic
$$dX_t=Y_t X_t dW_t,$$
where the volatility is given by 
$$dY_t= \mu (t, Y_t) dt + \sigma (t, Y_t)[\rho dW_t+\sqrt{1-\rho^2} dW^{'}_t].$$
Under these assumptions the different price sensitivities can be computed using subsections (3.1-3.2) and where the Malliavin derivative of $X_T$ and $Y_T$ are given by 
	For $0\leq t\leq T$, We have 
	\begin{equation*}
	D_t X_T=X_T \left(Y_t-\int_t^T Y_s D_t Y_s ds +\int_t^T D_t Y_s d{W}_s\right)
	\end{equation*}
	and
	\begin{equation*}
	D_t Y_s=\sigma (t, Y_t)\rho +\int_t^s \frac{\partial}{\partial y}(\rho\sigma+\mu)(\nu, Y_\nu) D_t Y_\nu d\nu,
	\end{equation*}
	here $\rho(\nu, Y_\nu)=\rho$.
For the $\alpha$-hypergeometric model of \cite{jose16}, the price of the underlying asset and its volatility are given by
$$dX_t=e^{Y_t} X_t dW_t,$$
where the volatility is given by 
$$dY_t= (a-be^{\alpha Y_t} )dt + \sigma dW^{'}_t,$$
where $dW_t dW^{'}_t:=\rho dt$. The Malliavin derivative of $X_T$ and $Y_T$ are
For $0\leq t\leq T$, We have 
\begin{equation*}
D_t X_T=X_T \left(e^{Y_t} -\int_t^T e^{2Y_s} D_t Y_s ds +\int_t^T e^{Y_s} D_t Y_s d{W}_s\right)
\end{equation*}
and
	\begin{equation*}
	D_t Y_s=\rho \sigma -\alpha b \int_t^s e^{Y_\nu}  D_t Y_\nu d\nu.
	\end{equation*}
The price sensitivities formulas obtained in (3.1) and (3.2) can be similarly applied to  most of the existing stochastic volatility models such as the Hull-White model \cite{hullwhite87} where the $Y$ process is lognormal, the correlation $\rho =0$ and the volatility is $f(y)=\sqrt{y}$, and the Heston model \cite{heston}, where $Y$ is a CIR process, the correlation $\rho\neq 0$ and the volatility is $f(y)=\sqrt{y}$.
\section{Conclusions}
This paper provides price sensitivities for a general stochastic volatility model that encompasses a number of recent models as well as well-known existing ones. For instance, the price sensitivities for recent stochastic volatility models such as the $\alpha$-hypergeometric model of \cite{jose16} and the linear model of \cite{jak13} can be found by applying the obtained results. These  results can be also used to find the price sensitivities for existing models that are used regularly by practitioners namely the Hull-White \cite{hullwhite87} and Heston \cite{heston} models. The Malliavin calculus is used for this purpose. The advantage of this method is that it is unbiased and it requires less computational time compare to the finite difference method which is commonly used in this context. 
%The suggested formulas are expressed as propositions combined with pertinent proofs. 
Each price sensitivity is a source of financial risk that investors need to tackle. Thus suggesting alternative measures of price sensitivities is expected to improve on the management of the underlying financial risk. 
%\section{Acknowledgement}
%The first author would like to express his thanks to the %College of Science at UAEU for supporting this research work %via the individual grant No. COS/IRG-12/15.
\section{Acknowledgment}
The first author would like to express his thanks to the College of Science at the UAE University for supporting this research work via the individual grant No. COS/IRG-12/15. The second author would like to acknowledge financial support by UAE University via UAEU-UPAR grant.

\end{document}